\documentclass[12pt]{article}

\usepackage{amsbsy,amsfonts,amsmath,amssymb,amsthm
,bbm,enumerate,euscript,graphicx}

\usepackage{amsthm}

\newtheorem{theo}{Theorem}
\newtheorem{coro}{Corollary}

\newtheorem{lemm}{Lemma}

\newcommand{\rems}{\paragraph{Remarks.}}

\newcommand{\lbl}{\label}

\newcommand{\eq}[1]{$(\ref{#1})$}

\newcommand{\Po}{{\cal P}}
\newcommand{\Q}{{\cal Q}}

\def\N{\mathbb{N}}
\def\Z{\mathbb{Z}}

\def\R{\mathbb{R}}

\def\Pr{\mathbb{P}}

\def\0{{\bf 0}}

\def\G{{\cal G}}
\def\cK{{\cal K}}
\def\cA{{\cal A}}
\def\cL{{\cal L}}

\newcommand{\toas}{\stackrel{{{\rm a.s.}}}{\longrightarrow}}

\newcommand{\eqco}{\setcounter{equation}{0}}
\newcommand{\thco}{\setcounter{theo}{0}}
\newcommand{\prco}{\setcounter{prop}{0}}
\newcommand{\laco}{\setcounter{lemm}{0}}
\newcommand{\coco}{\setcounter{coro}{0}}
\newcommand{\cjco}{\setcounter{conj}{0}}

\newcommand{\deco}{\setcounter{defn}{0}}
\newcommand{\allco}{\eqco  \thco \prco \laco \coco \cjco \deco}
\newcommand{\X}{{\cal X}}
\newcommand{\LL}{{\cal L}}

\newcommand{\A}{{\cal A}}

\newcommand{\Y}{{\cal Y}}
\newcommand{\I}{{\cal I}}

\renewcommand{\G}{{\cal G}}

\newcommand{\tl}{{\tilde{\lambda}}}

\newcommand{\eps}{\varepsilon}
\def\bdm{\begin{displaymath}}
\newcommand{\edm}{\end{displaymath}}
\def\benu{\begin{enumerate}}
\def\eenu{\end{enumerate}}
\def\beqn{\begin{equation}}
\def\eeqn{\end{equation}}
\def\be{\begin{equation}}
\def\ee{\end{equation}}
\def\bea{\begin{eqnarray}}
\def\eea{\end{eqnarray}}
\newcommand{\bean}{\begin{eqnarray*}}
\newcommand{\eean}{\end{eqnarray*}}
\newcommand{\bear}{\begin{eqnarray}}
\newcommand{\eear}{\end{eqnarray}}

\renewcommand{\epsilon}{\varepsilon}

\def\R{\mathbb{R}}

\renewcommand{\P}{{\mathbb P}}

\def\A{{\cal A}}

\def\qed{\hfill\hbox{${\vcenter{\vbox{
    \hrule height 0.4pt\hbox{\vrule width 0.4pt height 6pt
    \kern5pt\vrule width 0.4pt}\hrule height 0.4pt}}}$}}

\def\la{{\lambda}}

\begin{document}
\title{\bf Continuum AB percolation and AB random geometric
graphs}

\author{Mathew D. Penrose$^{1}$\\
{\normalsize{\em University of Bath}} }

\maketitle


 \footnotetext{ $~^1$ Department of
Mathematical Sciences, University of Bath, Bath BA2 7AY, United
Kingdom: {\texttt m.d.penrose@bath.ac.uk} }




\begin{abstract}
Consider a bipartite random geometric graph on
the union of two independent homogeneous Poisson point processes
in $d$-space, with distance parameter $r$ and intensities
$\lambda,\mu$.  We show for $d \geq 2$ that if $\lambda$ is  supercritical for
 the one-type random geometric graph with distance parameter $2r$,
 there exists
$\mu$ such that $(\lambda,\mu)$ is supercritical 
(this was previously known for $d=2$).
For $d=2$ we also consider the restriction of this graph to
points in the unit square. Taking $\mu = \tau \la$ for fixed $\tau$,
 we give a strong law of large numbers as $\la \to \infty$,
for the connectivity threshold of this graph.
\end{abstract}


\section{Introduction and statement of results}
The continuum AB percolation model, introduced by
Iyer and Yogeshwaran \cite{IY}, goes  as follows.
Particles of two types A and B are scattered 
randomly in  Euclidean space as 
two independent Poisson processes, and edges
are added between particles of opposite type
that are sufficiently close together. This
provides a continuum analogue to lattice AB percolation
which is discussed in e.g. \cite{Grimm}. Motivation
for considering continuum AB percolation is discussed in detail in \cite{IY};
the main motivation comes from wireless communications networks
with two types of transmitter.

Another type of continuum
percolation model with two types of particle is
the {\em  secrecy random graph} \cite{HS} in which the type B particles
(representing eavesdroppers) inhibit percolation; 
 each type A particle may send a message to every other
type A particle lying closer than its nearest neighbour of
type B. See also \cite{PW}. Such models are not considered here
but are complementary to ours.

To describe continuum AB percolation more precisely, we make some
definitions. Let $d \in \N$.
Given any two locally finite sets $\X,\Y \subset \R^d$,
and given $r >0$, 
let $G(\X,\Y,r)$ be the 
bipartite graph with vertex sets $\X$ and $\Y$, 
 and
with an undirected edge $\{X,Y\}$ included
for each $X \in  \X$ and $Y \in \Y$ with 
 $\|X-Y\| \leq r$, where $\|\cdot\|$ is the
Euclidean norm in $\R^d$ (our parameter $r$ would be denoted
by $2r$ in the notation of \cite{IY}.) Also,
let $G(\X,r)$ be the graph with vertex set $\X$
and
with an undirected edge $\{X,X'\}$ included
for each $X,X' \in  \X$ 
 with $\|X-X'\| \leq r$.

For $\lambda,\mu>0$  let $\Po_\lambda$, $\Q_\mu$ be
independent homogeneous Poisson point
 processes in $\R^d$ of intensity $\lambda, \mu $
respectively,  where
 we view each point process as a random subset of $\R^d$.
Our first results are concerned with the 
 bipartite graph 
 $G(\Po_\lambda , \Q_\mu,r)$.


Let $\I$ be the class of   graphs
having at least one infinite component.
By a version of the
Kolmogorov zero-one law,
given parameters $r,\lambda,\mu$ (and $d$), we have
 $\P[ G(\Po_\la, \Q_\mu,r ) \in \I] \in \{0,1\}$.
Provided $r, \lambda,$ and $\mu$ are sufficiently big, 
we have $\P[ G(\Po_\la, \Q_\mu,r ) \in \I] = 1$; see
\cite{IY}, or the discussion below.
Set
$$
\mu_c(r,\lambda): = \inf \{\mu:\P [ G(\Po_\lambda,\Q_\mu,r)  \in \I]
=1\},
$$
with the infimum of the empty set interpreted as $+\infty$.
Also, for  the more standard one-type continuum percolation
graph $G(\Po_\la,r)$,
define
$$
\lambda_c(2r): = \inf \{\lambda:\P [ G(\Po_\la,2r) \in \I]=1\},
$$
which is well known to be finite for $d \geq 2$ \cite{Grimm,MR},
 but is not  known analytically.
By scaling (see Proposition 2.11 of \cite{MR})
$\lambda_c(2r) = r^{-d} \lambda_c(2)$,
and explicit bounds for $\lambda_c(2)$ are provided in \cite{MR}.
Simulation studies indicate that $1- e^{-\pi \lambda_c(2)} \approx 0.67635$
for $d=2$ \cite{QZ} and
 $1 - e^{-(4\pi/3) \lambda_c(2)} \approx 0.28957$ for
$d=3$ \cite{LZ}.

Obviously if $G(\Po_\la,\Q_\mu,r)\in \I$ then also
$G(\Po_\la,2r) \in \I$,
 and therefore a necessary condition
for $\mu_c(r,\lambda)$ to be finite is that
$\lambda \geq \lambda_c(2r)$. In other words,
for any $r >0$ we have
\bea
\lbl{eq1}
\lambda_c^{AB}(r) :=
 \inf\{\lambda : \mu_c(r,\lambda) < \infty \} \geq \lambda_c(2r).
\eea
For  $d=2$ only, Iyer and Yogeshwaran \cite{IY} show that the inequality
in \eq{eq1} is in fact an equality.   For general $d \geq 2$,
 they also provide 
an explicit finite upper bound, here denoted $\tl_c^{AB}$,
  for $\lambda_c^{AB}(r)$, and
 establish
explicit upper bounds on $\mu_c(r,\lambda)$ for $\lambda > \tl_c^{AB}(r)$.
Note that even for $d =2$, their explicit upper bounds for 
 $\mu_c(r,\lambda)$ are given only when $\lambda > \tl_c^{AB}(r)$, with
 $\tl_c^{AB}(r) > \lambda_c(2r)$ for all $d \geq 2$; for 
the case with $d=2$ and $\lambda_c(2r) < \lambda \leq \tl^{AB}_c(r)$
their proof that $\mu_c(r,\lambda) < \infty$
 does not provide an explicit upper bound on
$\mu_c(r,\lambda) $.

In our first result, proved in Section \ref{secpf1}, we establish
for all dimensions (and all $r>0$)
 that the inequality in \eq{eq1}
is an equality, and provide explicit
asymptotic upper bounds  
on $\mu_c(r,\lambda)$ as  $\lambda $ approaches $ \lambda_c(2r)$
from above.
Let $\pi_d$ denote the volume of the unit radius ball in $d$ dimensions.
\begin{theo}
\lbl{th1}
Let $d \geq 2$ and $ r >0$. Then (i)
 $\lambda_c^{AB}(r) = \lambda_c(2r)$, and (ii)
writing $\lambda_c = \lambda_c(2r)$
we have
\bea
\lbl{eq2}
\limsup_{\delta \downarrow 0}
\left(
\frac{
\mu_c(r,\lambda_c + \delta) }{
\delta^{-2d} |\log \delta| } \right)
 \leq \left(\frac{4 \lambda_c^2}{r}
\right)^d d^{3d}(d+1) \pi_d.
\eea
\end{theo}
In  Section \ref{secpf1} we provide the proof, which is based on
the classic elementary continnum percolation techniques of discretisation,
 coupling and scaling. We shall also indicate how, for any
given $\lambda > \lambda_c(2r)$, one can
compute an explicit upper bound 
 for $\mu_c(r,\lambda)$; see eqn 
\eq{0918}.

It would be interesting to try to  find complementary {\em
lower} bounds for $\mu_c(r,\lambda)$.
An analogous problem in the lattice is
mixed  bond-site percolation, which similarly
has two parameters. For that model, similar questions
have been studied by Chayes and Schonman \cite{CS},
but it is not clear to what extent their methods
can be adapted to the continuum.

Our second result concerns full connectivity for
the {\em AB random geometric graph}, i.e. the
restriction of the AB percolation model to points
in a bounded region of $\R^d$.
 For $\la >0$ let $\Po_\la^F := \Po_\la \cap [0,1]^d$
and $\Q_\la^F := \Q_\la \cap [0,1]^d$
(these are finite Poisson processes of intensity $\lambda$;
hence the superscript $F$).
Given also $\tau >0$ and $r>0$,
let 
$\G^1(\la,\tau ,r)$
be the  graph on vertex set $\Po_\la^F $,
with an edge between each
pair of vertices sharing at least one common 
neighbour  in
$ G(\Po_\la^F , \Q_{\tau \la}^F ,r)$.

Let $\G^2(\la,\tau ,r)$
be the  graph on vertex set $\Q_{\tau \la}^F $,
with an edge between each
pair of vertices sharing at least one common 
neighbour  in
$ G(\Po_\la^F , \Q_{\tau \la}^F ,r)$.
Then
$ G(\Po_\la^F , \Q_{\tau \la}^F ,r)$
is connected,
if and only if both
 $\G^1(\la,\tau ,r)$ and
 $\G^2(\la,\tau ,r)$ are connected.

Let $\cK$ be the class of connected graphs, and let
\bean
\rho_n(\tau) = \min\{r: \G^1(n,\tau ,r) \in \cK \} 
\eean
which is a random variable determined by the
configuration of $(\Po_n,\Q_{\tau n})$. It is  a  {\em connectivity
threshold} for the AB random geometric graph.
Let us assume $\Po_\la^{F}$ and 
 $\Q_{\mu}^{F}$, are coupled for all $\la,\mu >0$ as follows. 
Let $(X_1,Y_1,X_2,Y_2,\ldots)$ be a sequence
of independent uniform random $d$-vectors uniformly
distributed over $[0,1]^d$. Independently, let $(N_t,t \geq 0)$ and
$(N'_t, t \geq 0)$
be independent Poisson counting processes of rate 1. Let
$\Po_\la^F = \{X_1,\ldots,X_{N_\la}\}$ and
$\Q_\mu^F = \{Y_1,\ldots,Y_{N'_\mu}\}$.

In Section \ref{seccon} we prove the following result, with
 $\toas$ denoting almost sure convergence  as $n \to \infty$ (with $n \in \N$).
\begin{theo}
\lbl{Thconn}
Assume $d=2$. Let $\tau >0$. Then
\bea
n \pi (\rho_n(\tau))^2 / \log n \toas \max (1/\tau,1/4). 
\lbl{0815e}
\eea
\end{theo}
\rems{
\begin{enumerate}
\item
The restriction to $d=2$ arises because boundary effects
become more important in higher dimensions (and $d=1$ is a different
case). It should be possible to adapt
the proof to obtain a similar result to \eq{0815e} 
in the unit {\em torus} in arbitrary dimensions $d \geq 2$,
namely $ n \pi_d (\rho_n(\tau))^d / \log n \toas \max (1/\tau,2^{-d}),$ 
although we have not checked the details.
\item
Iyer and Yogeshwaran \cite[Theorem 3.1]{IY} 
give almost sure lower and upper bounds for $\rho_n(\tau)$
in the torus. The extension of  our result
mentioned in the previous remark would
 show that the lower bound of \cite{IY}
 is sharp for $\tau \leq 2^d$,
and improve on their upper bound. 
\end{enumerate}
}

{\em Notation.} Given a countable set
 $\X $, we write $|\X|$ for the number of
elements of $\X$  and if also
 $\X \subset \R^d$, given
$A \subset \R^d$ we write $\X(A)$ for $|\X \cap A|$. 
Also, for  $a>0$ we write $aA$ for
$\{ ay: y \in A\}$.
Let $\oplus$ denote Minkowski addition of sets (see
 e.g. \cite{Penbk}).

\section{Percolation: proof of Theorem \ref{th1}}
\lbl{secpf1}
\allco
%
Fix $r>0$ and 
let $\lambda > \lambda_c(2r)$.
We first prove that $\mu_c(r,\lambda) < \infty$; combined with \eq{eq1}
this shows that 
 that $\lambda_c^{AB}(r) = \lambda_c(2r)$, which is part (i) of the theorem.
Later we shall quantify the
estimates in our argument, thereby establishing  part (ii).

Choose $s <r$ and $\nu < \lambda$ such that $\P[G(\Po_\nu,2s) \in
\I ]=1$.
This is possible because 
 decreasing the radius slightly is equivalent to decreasing
 the Poisson intensity slightly, by scaling
(see \cite{MR}; also
 the first equality of \eq{scaleq} below). Set $t = (r+s)/2$, and
let $\eps >0$ be chosen small enough
so that any cube of side $\eps$  has Euclidean diameter at most $t-s$. 
For $a >0$ 
let $p_a:= 1 - \exp(- \eps^d a)$, the probability that  a given cube of side
$\eps$ contains at least one point of $\Po_a$.

Consider Bernoulli site percolation on the graph $(\eps \Z^d, \sim)$
where for $u,v \in \eps \Z^d$ we put $u \sim v$ if and only if there exists
$w \in \eps \Z^d$ with $\|w-u\| \leq t$ and  $\|w-v\| \leq t$.
Given $p >0$
suppose each site $u \in \eps \Z^d$ is
 independently occupied with probability $p$.
Let $D_1$ be the event that there is an infinite
path of occupied sites in the graph, and let $\P_p[D_1]$
be the probability that this event occurs.

 Divide $\R^d$ into cubes $Q_u, u \in \eps \Z^d$,
defined by $Q_u := \{u\} \oplus [0,\eps)^d$.
For $x \in \R^d$ let $z_x \in \eps \Z^d $ be such that 
$x \in Q_{z_{x}}$.
The Poisson 
 process
 $\Po_\nu$
 may be coupled to
 a realization of
the site percolation process with parameter $p_\nu$,
 by deeming each $z \in  \eps \Z^d$
to be occupied  if and only if $\Po_\nu(Q_{z}) \geq 1$.
By the choice of $\eps$, for $X,Y \in \Po_\nu$, if $\|X-Y\| \leq 2s$ then
$\|z_X - z_{(X+Y)/2}\| \leq t$ and 
$\|z_Y - z_{(X+Y)/2}\| \leq t$, and hence $z_X \sim z_Y$.
Therefore, with this coupling, if  $G(\Po_\nu,2s) \in \I$
then there is an infinite path of occupied sites in $(\eps\Z^d,\sim)$. Since
 we chose $\nu, s$ so that $\P[G(\Po_\nu,2s) \in \I] =1$,
 we have $\P_{p_\nu}[D_1] =1$.

Now consider a form of lattice AB percolation on $\eps \Z^d$
with parameter  pair $(p,q) \in [0,1]^2$ (not necessarily the
same as any of the lattice AB percolation models in the literature). 
Let $(V_u, u \in \eps \Z^d)$ be a family of independent Bernoulli
random variables with parameter $p$, and
let $(W_u, u \in \eps \Z^d)$ be a family of independent Bernoulli
random variables with parameter $q$. Let $D_2$
be the event
that there is an infinite sequence $u_1,u_2,\ldots,$ of
distinct elements of $\eps \Z^d$, and
 an infinite sequence $v_1,v_2,\ldots,$ of
 elements of $\eps \Z^d$, such that
for each $i \in \N$ we have $V_{u_i}W_{v_i}=1$ and
$\max(\|u_i-v_i\|,\|v_i-u_{i+1}\|)\leq t$.
Let $\tilde{\P}_{p,q}[D_2]$ be the probability that event $D_2$ occurs,
given the parameter pair $(p,q)$.


Since $\P_{p_\nu}[D_1]=1$, clearly
 $\tilde{\P}_{p_\nu,1}[D_2]=1$.
Increasing $p $ slightly and decreasing $q$ slightly,
we shall show that there exists $q<1$ such that
 \bea
\tilde{\P}_{p_\lambda,q}[D_2]=1.
\lbl{eq3}
\eea
This is enough to demonstrate that $\mu_c(r,\lambda)< \infty$.
Indeed, suppose such a $q$ exists and 
 choose $\mu$ such that $p_\mu =q$.
Then for $u \in \eps \Z^d$ set
 $V_u = 1 $ if and only if $\Po_\lambda(Q_u) \geq 1$ and 
$W_u = 1 $ if and only if $\Q_\mu(Q_u) \geq 1$.
Suppose 
$D_2$ occurs and let $u_1,v_1,u_2,v_2,\ldots$ be
as in the definition of event $D_2$.
Then for each $i\in \N$ we have
$V_{u_i}=1$ so we can pick a point $X_i \in \Po_\lambda \cap Q_{u_i}$,
and
$W_{v_i}=1$ so we can pick a point $Y_i \in \Q_\mu \cap Q_{v_i}$.
Then by the choice of $\eps$, for each $i \in \N$ we have
$$
\max(\|X_i-Y_i\|,\|Y_i-X_{i+1}\|) \leq t+ (t-s) = r,
$$
and hence $G(\Po_\la, \Q_\mu,r) \in \I$.
  Hence, by \eq{eq3} we have $\P [G(\Po_\lambda,\Q_\mu,r) \in \I ]=1$.
Therefore $\mu_c(r,\lambda ) \leq \mu < \infty$ as asserted.

To complete the proof of part (i), it remains
 to prove that  \eq{eq3} holds for some $q <1$.
Let  $(T_u,u \in \eps \Z^d)$
be independent Bernoulli variables with parameter $p_\lambda$.
For each ordered pair $(u,v) \in (\eps \Z^d)^2$ with $0<\|u-v\| \leq t$, let
$U_{u,v}$ be independent Bernoulli random variables with
parameter $(p_\nu/p_\lambda)^{1/\Delta}$,
where we set
 \bea
\Delta := | \{u \in \eps \Z^d: 0< \|u\| \leq t \} |.
\lbl{Deltadef}
\eea
 Assume   the variables $U_{u,v}$ and $T_u$
are all mutually independent.
Then for $u,v \in \eps \Z^d$  define the Bernoulli variables
\bea
V_u : = T_u \prod_{v \in \eps \Z^d: 0 < \|v-u\| \leq t } U_{u,v} ;
\\
W_v:= 1-
 \prod_{u \in \eps \Z^d: 0 < \|v-u\| \leq t }
 \left(1-
 U_{u,v} \right).
\eea
Then $(V_u)_{u \in \eps \Z^d}$ are a family of
 independent Bernoulli variables with parameter $p_\nu$. 
Also,
 $(W_v)_{v\in \eps \Z^d}$
 are a family of independent Bernoulli
variables with parameter
 \bea
q : = 1 -(1-(p_\nu/p_\lambda)^{1/\Delta})^\Delta < 1,
\lbl{eq4}
\eea
and are independent of $(T_u, u \in \eps \Z^d)$.

Since $\P_{p_\nu}[D_1] =1$,
 with probability 1 there exists an infinite
sequence $u_1,u_2,\ldots $ of distinct elements
of $\eps \Z^d$ with $u_i \sim u_{i+1}$ for all $i \in \N$,
and with $V_{u_i} =1$ for each $i \in \N$.
By definition of the relation $\sim$, we can
choose  sequence $v_1,v_2,\ldots$ of elements of $\eps \Z^d$
such that for each $i \in \N$ we have $\max(\|v_i -u_i\|,\|v_i-u_{i+1}\|) \leq t$.
Then for each $i$, since $V_{u_i} =1$ we have $U_{u_i, v_i} =1$, and therefore
$W_{v_i}=1$; also $T_{u_i} =1$.
 Hence, \eq{eq3} holds
as required,
establishing 
 that $\mu_c(r,\lambda) < \infty$. This completes the proof of
part (i). \\

For part (ii), we need to quantify the preceding argument. 
First note that the value of $\mu$ associated
with $q$ given by \eq{eq4} (i.e. with $p_\mu =q$)
has 
$
\exp(-\mu \eps^d) = 
(1-(p_\nu/p_\lambda)^{1/\Delta})^\Delta,
$
so that since $ \eps^d \Delta \leq \pi_d r^d$ by \eq{Deltadef}, we have 
\bea
\mu_c (r,\lambda) \leq \mu 
= \eps^{-d} \Delta \log \left( \frac{1}{1 - 
(p_\nu/p_\lambda)^{1/\Delta} } \right) 
\nonumber \\
 \leq
 \eps^{-2d} \pi_d r^d \log \left( \frac{1}{1 - 
(p_\nu/p_\lambda)^{(\eps/r)^d/\pi_d} } \right). 
\lbl{eq5}
\eea
From now on set $\lambda_c  := \lambda_c(2r)$, 
and set $\lambda = \lambda_c + \delta$ for some $\delta >0$.
We need to choose $s <r$ and $\nu < \lambda$ such
that 
 $\P[G(\Po_\nu,2s) \in \I] =1$.
Choose $\alpha, \beta >0$ with $\alpha + \beta <1$,
and also let $\alpha' \in (0,\alpha)$ and $\beta' \in (0,\beta)$.
Set
$$
s: =
r ( 1 + \alpha \delta/\lambda_c)^{-1/d}; ~~~~~
\nu := \lambda_c + (1 - \beta) \delta.
$$
By scaling 
(see Proposition 2.11 of \cite{MR}) 
and our choice of $s$, we have 
\bea
\lambda_c(2s) = (r/s)^d \lambda_c(2r)
= \lambda_c + \alpha \delta,
\lbl{scaleq}
\eea
and hence $\nu > \lambda_c(2s)$ so 
 $\P[G(\Po_\nu,2s) \in \I] =1$,
as required.

 Our choice of $\eps$ in the discretization needs to satisfy
\bea
\eps \leq \frac{r-s}{2 \sqrt{d}} 
= \frac{r}{2 \sqrt{d}}
 \left[ 1 - \left(
1+
 \frac{\alpha \delta}{ \lambda_c}
\right)^{-1/d}
 \right],
\lbl{eq6}
\eea
and the right hand side of \eq{eq6} is asymptotic to
$\alpha r \delta /(2 d^{3/2} \lambda_c) $ 
as
 $\delta \to 0$. Hence, taking
$\eps = \alpha' r\delta/(2 d^{3/2} \lambda_c)$, we have \eq{eq6}
provided $\delta \leq \delta_1$, for some fixed $\delta_1 >0$.
Also,
\bea
 \frac{p_\nu}{p_\lambda} 
\leq 
 \frac{ \eps^d \nu}{ \eps^d \lambda \exp(-\eps^d \lambda)}
=
 \left(\frac{\lambda_c + (1-\beta) \delta }{\lambda_c +  \delta}\right)
\exp(\eps^d \lambda) 
\lbl{0908a}
\eea
 and so by Taylor expansion, there is some $\delta_2 >0$ such that
provided $0< \delta \leq \delta_2$, taking
$\eps = \alpha' r\delta/(2 d^{3/2} \lambda_c)$
 we have
\bea
\left( \frac{p_\nu}{p_\lambda} \right)^{(\eps/r)^d/\pi_d}
\leq 1 - \beta' \delta \eps^d /(  \pi_d r^d \la_c)
= 1 - \frac{\beta' \delta^{d+1}}{ \pi_d  \la_c(2 d^{3/2} \lambda_c/\alpha')^d}.
\nonumber
\eea
Therefore by \eq{eq5}, for $0 < \delta \leq  \min(\delta_1, 
\delta_2)$ we have
\bean
\mu_c(r,\lambda) \leq \left(  \frac{2 d^{3/2}\lambda_c}{r \delta \alpha'}
 \right)^{2d}
\pi_d r^d \log[ \pi_d  \la_c(2 d^{3/2} \lambda_c/\alpha')^d (1/\beta') \delta^{-d -1} ]
\eean
and since we can take $\alpha'$ arbitrarily close to 1,
 \eq{eq2} follows, completing the proof. \\


For a given value of $\lambda$ with $\lambda = \lambda_c(2r) + \delta$
for some $\delta >0$, an explicit upper bound for
$\mu_c(r,\lambda)$ could be computed as follows. Choose  $\alpha, \beta
> 0$ with $\alpha + \beta <1$, and let $\eps$ be given by the right hand side of
\eq{eq6}. 
%
%
%
%
Then a numerical upper bound for $\mu_c(r,\lambda)$ can
be obtained by computing the right hand side of \eq{eq5}.
To make this bound as small as possible (given $\alpha$), we make
$\nu$ as small as we can, i.e. 
 make $\beta$ approach  $ 1-\alpha$ and $\nu $ approach $\lambda_c
+  \alpha\delta$.
Taking this limit and then optimizing further over $\alpha$ gives us
the upper bound
\bea
\mu_c (r,\lambda) \leq  \inf_{\alpha \in (0,1)}
 \eps(\alpha)^{-2d} \pi_d r^d \log \left( \frac{1}{1 - 
(p_{\lambda_c +\alpha \delta}/p_\lambda)^{(\eps(\alpha)/r)^d/\pi_d} } \right), 
\lbl{0918}
\eea
with $\eps = \eps(\alpha)$ given by the right side of \eq{eq6}.

\section{Connectivity: proof of Theorem \ref{Thconn}}
\lbl{seccon}

\allco

Throughout this section we assume $d=2$.
 All asymptotics
are as $n \to \infty$. Given $a, b \in \R$
we shall sometimes write $a \vee b$ for $\max(a,b)$ and
 $a\wedge b$ for $\min(a,b)$. 
Fix $\tau >0$.
Given $\tau$ and $r_n$, let $\delta_n$ denote the 
minimum degree of $\G^1(n,\tau,r_n)$. 
\begin{lemm} \lbl{Thmin}
Let $\alpha \in ( 0,1/\tau)$. If
 $n \pi r_n^2/\log n = \alpha $, $n \geq 2$,
then almost surely, $\delta_n = 0$ for all
but finitely many $n$.
\end{lemm}
{\em Proof.}
See \cite[Proposition 5.1]{IY}.
\begin{lemm} \lbl{Thmin2}
Let $\alpha \in (0,1/4)$.
If $n \pi  r_n^2/\log n = \alpha, n \geq 2$,
then almost surely, $\delta_n = 0$ for all
but finitely many $n$.
\end{lemm}
{\em Proof.}
By  \cite[Theorem 7.8]{Penbk}, for this choice of $r_n$,
almost surely
 the minimum degree of the (one-type) 
geometric  graph $G(\Po_n^F ,2r_n)$
is zero for all but finitely many $n$,
 and therefore so is the minimum degree of
 $\G^1(n,\tau,r_n)$. $\qed$

\begin{coro}
Let $d=2$.
Given $\eps >0$ we have almost surely that
 $
 n \pi (\rho_n(\tau))^2/ \log n > (1- \eps) \max(1/4,1/\tau) 
$ for all but finitely many $n$.
\end{coro}
{\em Proof.}
Assume $\eps <1$.  
For $n \geq 2$, set $r_n = ((1-\eps) (1/4 \vee 1/\tau) \log n/(n\pi))^{1/2}$,
so $n \pi r_n^2/\log n = (1- \eps)(1/4 \vee 1/\tau)$.
Let $\delta_n$ be the minimum degree of $\G^1(n,\tau,r_n)$.
If the minimum degree of a graph of order greater than 1 is zero,
then it is not connected; hence
\bean
 \{ n \pi (\rho_n(\tau))^2/ \log n \leq ( 1- \eps ) ( \frac{1}{4} \vee
\frac{1}{\tau} ) \} 
= \{ \G^1( n,\tau,r_n) \in \cK \}
\\
\subset \{ \delta_n >0 \} \cup \{ \Po_n^F([0,1]^2 ) \leq 1 \},
\eean
which occurs only finitely often almost surely,
 by Lemmas \ref{Thmin} and \ref{Thmin2}. $\qed$ \\

To complete the proof of Theorem \ref{Thconn},
it suffices to prove the following:
\begin{theo}
\lbl{connth2}
Suppose for some fixed $\alpha $ that $(r_n)_{n \in \N}$ is such that
for all $n \geq 2$,
\bea
n \pi r_n^2/\log n = \alpha > \max (1/\tau,1/4).
\lbl{0818a}
\eea
Then almost surely  $\G^1(n,\tau,r_n) \in \cK$ for all
but finitely many $n$.
\end{theo}
The proof of this requires a series
of lemmas. It proceeds by discretization of
space.
Assume $\alpha$ and $r_n$ are given, satisfying \eq{0818a}.
Let $\eps_0 \in(0,1/99)$ be chosen in such a way that
for $\eps = \eps_0$ we have
\bea
\alpha \tau ( 1- 12 \eps ) > 1+ \eps;
\lbl{epscond}
\\
\alpha (4 - 12 \eps (3 + \tau) ) > 1+ \eps.
\lbl{epscond2}
\eea
Given $n$,
partition $[0,1]^2$ into squares of side $\eps_n r_n$ with $\eps_n  $
chosen so that $\eps_0 \leq \eps_n < 1/99$ and
 $1/(\eps_n r_n) \in \N$, and $\eps = \eps_n$ satisfies 
\eq{epscond} and \eq{epscond2};
this is possible for all large enough $n$, say for $n \geq n_0$. 
In the sequel  we assume $n \geq n_0$ and
 often write just $\eps $ for $\eps_n$.

Let $\cL_n$ be the set of centres of the squares in this
partition (a finite lattice).
Then $|\cL_n| = \Theta ( n/\log n)$.
List the squares as $Q_i, 1 \leq i \leq |\cL_n|$,
and the corresponding centres of squares (i.e., the elements
of $\cL_n$) as $q_i, 1 \leq i \leq |\cL_n|$.

Given a set $\X \subset [0,1]^2$, define
the {\em projection of $\X$ onto $\LL_n$}
to be the set of $q_i \in \LL_n$  such 
that $\X \cap Q_i \neq \emptyset$.
Given also $\Y \subset [0,1]^2$, 
 define the projection of $(\X,\Y)$ onto $\LL_n$ to be the pair
 $(\X',\Y')$, where
 $\X'$ is the projection of $\X$ onto $\LL_n$ and 
 $\Y'$ is the projection of $\Y$ onto $\LL_n$.
We refer to 
$|\X'|+|\Y'|$ (respectively 
$|\X'|$, $|\Y'|$) as the
 {\em order} of  the projection of $(\X,\Y)$
(respectively of $\X$, of $\Y$) onto $\LL_n$.

\begin{lemm}
\lbl{conlem1}
Let $n \in \N$.
Suppose $\X$ and $\Y$ are
finite subsets of $[0,1]^2$, such
that $G(\X,\Y,r_n) $ is connected.
Let $(\X',\Y')$ be the projection of $(\X,\Y)$
onto $\LL_n$.
Then the bipartite geometric graph $G(\X',\Y',r_n(1+2 \eps_n))$ is
connected.
\end{lemm}
{\em Proof.}
If $q_i,q_j \in \cL_n$
and  $X \in \X$ and $Y \in \Y$,
 with $\|X -Y\| \leq r_n$,
then 
by the triangle inequality
we have 
$$
\|q_i-q_j\| \leq \|X-q_i\| + \|X -Y \| + \|Y-q_j\| \leq 
 r_n(1+2 \eps)
$$
and therefore since $G(\X,\Y,r_n)$ is connected,
so is $G(\X',\Y',r_n(1+2 \eps))$.
$\qed$\\

Given $n,m \in \N$, let $\cA_{n,m}$ denote the set of pairs
 $(\sigma_1,\sigma_2)$ with each $\sigma_j \subset \cL_n$,
with $|\sigma_1|+|\sigma_2| =m$ and $|\sigma_1| \geq 1$,
 such that $G(\sigma_1,\sigma_2,r_n(1+2\eps_n))$
is connected; these may be viewed as `bipartite lattice animals'.

Let $\cA_{n,m}^2 $ be the set of $(\sigma_1,\sigma_2) \in \cA_{n,m}$
such that
 all elements of $\sigma_1 \cup \sigma_2$ are distant at least $2r_n$ from
the boundary of $[0,1]^2$.

Let $\cA_{n,m}^1 $ be the set of $(\sigma_1,\sigma_2) \in \cA_{n,m}$
such that
 $\sigma_1 \cup \sigma_2$ is distant less than  $2r_n$ from
{\em just one edge} of $[0,1]^2$.

Let $\cA_{n,m}^0 :=
\cA_{n,m} \setminus ( \A_{n,m}^2 \cup \A_{n,m}^1 )$,
the set of $(\sigma_1,\sigma_2) \in \cA_{n,m}$
such that
 $\sigma_1 \cup \sigma_2$ is distant less than  $2r_n$ from
{\em two  edges} of $[0,1]^2$ (i.e. near a corner of $[0,1]^2$).

\begin{lemm}
\lbl{countlem}
Given $m \in \N$, there is constant $C = C(m)$ such that
for all $n \geq n_0$ we have
$$
|\cA_{n,m} | \leq C(n/\log n), ~~~~
|\cA_{n,m}^1  |\leq C(n/\log n)^{1/2}, ~~~~
|\cA_{n,m}^0 | \leq C.
$$
\end{lemm}
{\em Proof.}
Fix $m$.
Consider how many ways there are to choose
$\sigma \in \cA_{n,m}$.

There are at most $|\LL_n|$ choices, and hence $O(n/\log n)$ choices,
 for the first
element of $\sigma_1$ in the lexicographic ordering.
Having chosen the first element of $\sigma_1$, there are a bounded
number of ways to choose the rest of $\sigma$.

Consider how many ways there are to choose
$\sigma \in \cA_{n,m}^1$.
In this case there are $O(r_n^{-1} ) = O(( n / \log n)^{1/2})$
ways to choose the first element of $\sigma_1$ (distant at most
$2 r_n$ from the boundary of $[0,1]^2$), and then
a bounded number of ways to choose the rest of $ \sigma$.

Finally consider how many ways there are to choose
$\sigma \in \cA_{n,m}^0$.
In this case there are $O(1) $
ways to choose the first element of $\sigma_1$, and then
a bounded number of ways to choose the rest of $ \sigma$.
$\qed$ \\

For $n \in \N$ set $\nu(n) := n^{\lceil 4/\eps_0 \rceil}$.
Note that $\nu(n+1) \sim \nu(n)$ and
$r_{\nu(n+1)} \sim r_{\nu(n)}$ as $n \to \infty$,
and that
$r_n$ is monotone decreasing in $n$, $n \geq 3$.

Given $n \in \N$ with $\nu(n) \geq n_0$, 
and given $\sigma_1 \subset\cL_{\nu(n)}$
and $\sigma_2 \subset \cL_{\nu(n)}$,
 let $E_{(\sigma_1,\sigma_2)}$
be the event that  there exists 
 some $n' \in \N \cap [\nu(n),\nu(n+1))$
such that there is
 a component $(U,V)$ 
 of $G(\Po_{n'}^F,\Q_{\tau n'}^F, r_{n'})$
 such
that $(\sigma_1,\sigma_2)$ is the projection of $(U,V)$ onto 
$\LL_{\nu(n)}$.

For $x \in \R^2$ and $r >0$ let 
$B(x,r) := \{y \in \R^2:\|y-x\| \leq r\}$. Also
let $B_+(r) $ be the right half of $B((0,0),r)$, 
and let $B_-(r) $ be the left half of $B((0,0),r)$. 
Let $v_2(\cdot)$ denote Lebesgue measure,
defined on Borel subsets of $\R^2$.

\begin{lemm}
\lbl{AniLem}
There exists $n_1 \in \N$ such that
 for all $m \in \N$ and
 $n \geq n_1$ we have
\bea
\sup_{\sigma \in \A_{\nu(n),m}^2} ( \Pr[E_\sigma] )  \leq \nu(n)^{-(1+\eps)}; 
\lbl{Ani1}
\\
\sup_{\sigma \in \cA_{\nu(n),m}^1} (  \Pr[E_\sigma] ) \leq \nu(n)^{-(1+\eps)/2}; 
\lbl{Ani2}
\\
\sup_{\sigma \in \cA_{\nu(n),m}^0} (\Pr[E_\sigma]) \leq \nu(n)^{-1/20}.  
\lbl{Ani3}
\eea
\end{lemm}
{\em Proof.}
Choose $n_1$ so that $\nu(n_1) \geq n_0$ and also
$(1-\eps_0) r_{\nu(n)}  < r_{\nu(n+1)}$ for $n \geq n_1$.
Assume from now on that $n \geq n_1$.

Given $\sigma =(\sigma_1,\sigma_2) \in \cA_{\nu(n),m}$, let $q_i$
(respectively $q_j$) be the lexicographically
first (resp. last) element of $\sigma_1$. 
Let $\sigma_2^-$ be the set of 
$q_k \in \sigma_2 \cap B(q_i,r_{\nu(n)}(1-4 \eps))$
lying strictly to the left of $q_i$ (in this proof, $\eps := \eps_{\nu(n)}$).
Let $\sigma_2^+$ be the set of 
$q_k \in \sigma_2 \cap B(q_j,r_{\nu(n)}(1-4 \eps))$
lying strictly to the right  of $q_j$.
Let $\tilde{\sigma}_2^+ : = \sigma_2^+ \oplus [-\eps r_{\nu(n)}/2,\eps
 r_{\nu(n)}/2]^2$ and
 $\tilde{\sigma}_2^- : = \sigma_2^- \oplus [-\eps r_{\nu(n)}/2,\eps 
r_{\nu(n)}/2]^2$ (see Figure 1).

Let $B_{\sigma}^{-}$ be the part of $B(q_i,r_{\nu(n)}(1- 5 \eps))$ lying 
strictly to the left of
$Q_i$. Let $B_{\sigma}^{+}$ 
be the part of $B(q_j,r_{\nu(n)}(1- 5 \eps))$ lying 
strictly to the right of
$Q_j$. 

Given $\sigma$,
define the events $A_\sigma^+$ 
and  $A_\sigma^-$ 
by
\bean
A_\sigma^+ :=
\{ \Q_{\tau \nu(n+1)}^F(B_\sigma^+ \setminus \tilde{\sigma}_2^+) = 0\}
 \cap
 \{ \Po_{\nu(n+1)}^F(\tilde{\sigma}_2^+ \oplus B_+(r_{\nu(n)}(1 - 3 \eps) )) = 0 \};
\\
A_\sigma^- :=
\{ \Q_{\tau \nu(n+1)}^F(B_\sigma^- \setminus \tilde{\sigma}_2^-) = 0\}
 \cap
 \{ \Po_{\nu(n+1)}^F(\tilde{\sigma}_2^- \oplus
 B_-(r_{\nu(n)}(1 - 3 \eps) )) = 0 \}.
\eean
See Figure 1 for
an illustration of event $A_\sigma^+$.
Note that events $A_\sigma^+$ and $A_\sigma^-$ are independent.

\begin{figure}[!h]
\center
\includegraphics[width=0.5\textwidth]{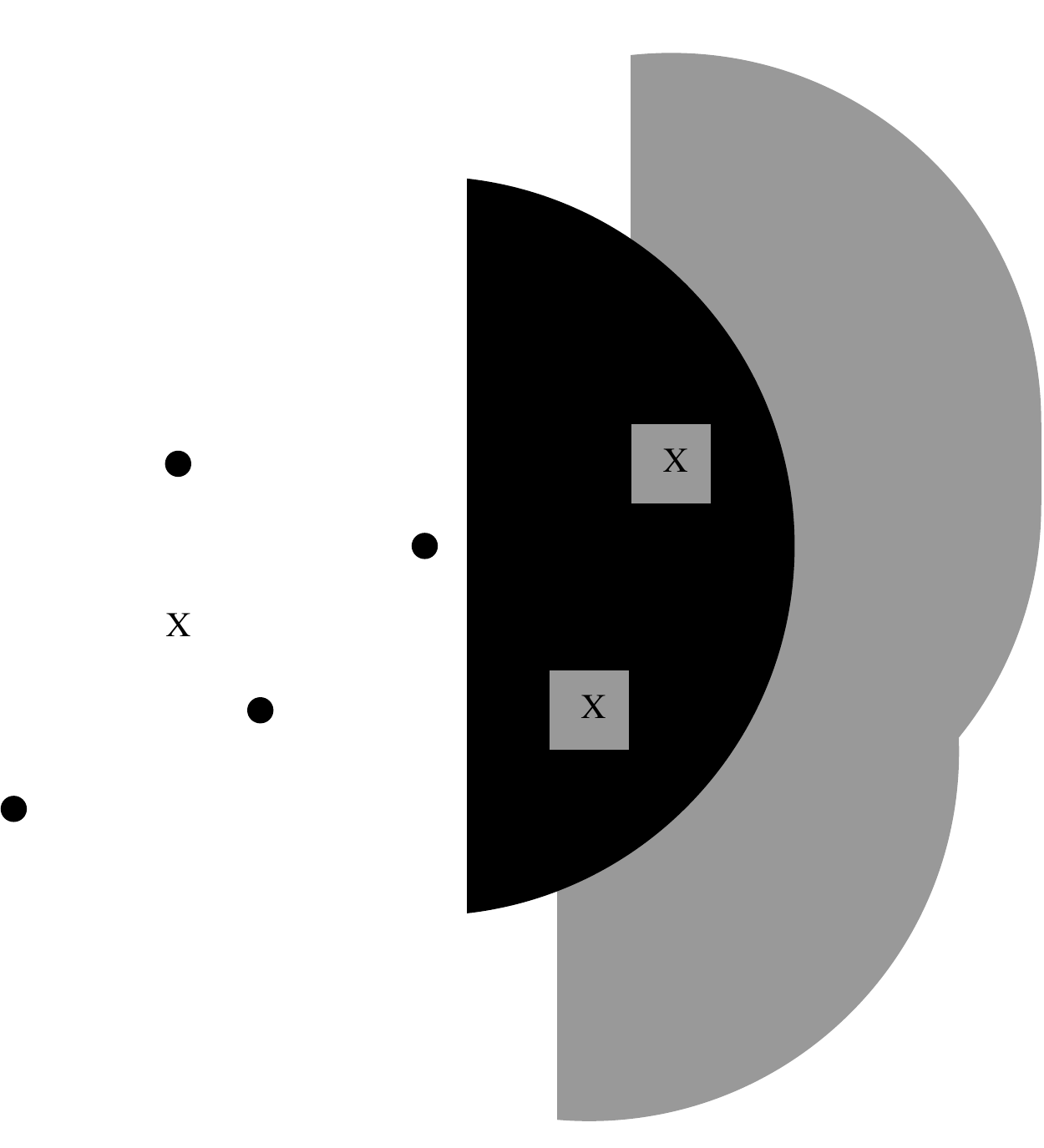}
\caption{
The dots are the points of $\sigma_1$, and
the crosses are the points of $\sigma_2$.
The grey squares are the set $\tilde{\sigma}_2^+$
(since $\eps = \eps_n < 1/99$, they should really be
smaller).
Event $A_\sigma^+$ says that
 the black region contains no points of $\Q_{\tau \nu(n+1)}^F$ and
 the grey region (partly obscured by the black region)
 contains no points of $\Po_{\nu(n+1)}^F$.
}
\end{figure}


Suppose $k$ is such that 
$Q_k \cap  B_\sigma^+ \neq \emptyset$.
Then  by the triangle equality,
\bea
\|q_k-q_j\| \leq r_{\nu(n)}(1-5 \eps) + \eps r_{\nu(n)} = r_{\nu(n)}(1- 4 \eps).
\lbl{0903a}
\eea 
Similarly, if 
$Q_k \cap  B_\sigma^- \neq \emptyset$
then $\|q_k-q_i\| \leq r_{\nu(n)}(1-4 \eps)$.

By our coupling of Poisson processes,
 for $\nu(n) \leq n' < \nu(n+1)$ we have
 $\Po_{\nu(n)} \subset \Po_{n'} \subset \Po_{\nu(n+1)}$.
Also if $x \in Q_k$ and $y \in Q_i$ with $\|q_i -q_k\| \leq r_{\nu(n)}(1- 3 \eps)$,
then by the triangle inequality and our condition on $n_1$
we have $\|x-y\| \leq r_{\nu(n)}(1- \eps) \leq r_{n'}$ for
all $n' \in  [\nu(n),\nu(n+1))$.
Hence by the argument at \eq{0903a}, for any $\sigma \in \A_{n,m}$
we have
$E_\sigma \subset A_\sigma^+ \cap A_\sigma^-$.

First  we prove \eq{Ani2}. Take $\sigma \in \A_{\nu(n),m}^1$.
Consider just the case where $\sigma $ is near to the left 
edge of $[0,1]^2$ (the other three cases are 
treated similarly).
If $\sigma_2^+ = \emptyset$, then
$A_\sigma^+ = \{ \Q_{\tau \nu(n+1)}^F (B^+_\sigma)=0 \}$, 
and in this case
 we have
\bea
\Pr[A_\sigma^+] 
\leq \exp ( - \tau \nu(n) (\pi (r_{\nu(n)}(1-5 \eps))^2 
- 2 \eps r_{\nu(n)}^2) /2 )
\nonumber \\
\leq \exp \left[ -\tau \alpha (\log \nu(n)) 
( 1- 12 \eps) /2 \right] \leq  
 \nu(n)^{-(1+\eps)/2},
\lbl{0903c}
\eea
where the last inequality comes from
  \eq{epscond}.  This proves  \eq{Ani2} for this case.

Suppose instead that $\sigma_2^+  \neq \emptyset$. 
Then $\tilde{\sigma}_2^+ \subset \{q_j\} \oplus B_+(r_{\nu(n)}(1- 3 \eps))$
so that $v_2(\tilde{\sigma}_2^+ ) \leq \pi r_{\nu(n)}^2 (1- 3 \eps)^2/2$.
Let
 $s \in [0,1]$ be chosen such that
$v_2(\tilde{\sigma}_2^+) = s^2 \pi r_{\nu(n)}^2 (1- 3 \eps)^2/2$.
Then by the Brunn-Minkowski inequality (see e.g. \cite{Penbk}),
$$
 v_2( \tilde{\sigma}_2^+ \oplus B_+(r_{\nu(n)}(1- 3 \eps))) \geq
(\pi r_{\nu(n)}^2/2)(1-  3 \eps)^2 ( 1+ s)^2, 
$$
and also $v_2(B_\sigma^+) \geq \pi r_{\nu(n)}^2((1- 5 \eps)^2 - 2 \eps)/2$
so that
\bean
\Pr[A_\sigma^+]
\leq \exp ( - \tau \nu(n) v_2( B_\sigma^+ \setminus \tilde{\sigma}_2^+)
- \nu(n) v_2( \tilde{\sigma}_2^+ \oplus B_+(r_{\nu(n)}(1- 3 \eps))))
\\
\leq \exp( - (\nu(n) \pi r_{\nu(n)}^2/2)[
 \tau ((1- 5 \eps)^2 -2 \eps -s^2(1-3 \eps)^2)  + (1 + s)^2(1-3 \eps)^2 ] )
\\
\leq 
\exp (- (\alpha/2) (\log \nu(n)) g_\tau (s)
),
\eean
where we set $g_\tau(s) := ( \tau + 1 + 2s)(1- 12 \eps) 
+ s^2(1- 3 \eps)^2 (1-\tau)$.
If  $\tau \leq 1$ then $g_\tau(s)$  
is
minimised over $s \in [0,1]$ at $s =0$. 
If
 $\tau > 1$, then $g_\tau(\cdot)$
is concave, so its minimum over 
$[0 , 1]$ is achieved at $s =0 $ or $s=1$;
also in this case $g_\tau(1) \geq (3 + \tau) (1- 12 \eps) + 1-\tau$.
Hence,
 using \eq{epscond} and \eq{epscond2} we obtain
\bea
\Pr [ A_\sigma^+ ] \leq \exp (- (\alpha/2)  (\log \nu(n))
\min [ (1+ \tau) (1 -12 \eps) ,4 - 12 \eps (3+\tau)] )
\nonumber \\
\leq \nu(n)^{-(1+\eps)/2},
\lbl{0903b}
\eea
completing the proof of \eq{Ani2}.

Now we prove \eq{Ani1}.
If $\sigma \in \A^2_{n,m}$ then 
$\Pr [ A_\sigma^+ ] \leq \nu(n)^{-(1+\eps)/2}$ by
\eq{0903c} and
\eq{0903b}, and
$\Pr [ A_\sigma^- ] \leq \nu(n)^{-(1+\eps)/2}$ similarly.
Therefore
$
\Pr[ E_\sigma ] \leq \Pr [ A_\sigma^+ \cap A_\sigma^-] \leq 
 \nu(n)^{-1 -\eps},
$
completing the proof of \eq{Ani1}.

 To prove \eq{Ani3}, let $\sigma \in \A_{n,m}^0$.
Assume $\sigma$ is near the lower left corner of
$[0,1]^2$ (the other cases are treated similarly).
First suppose $\sigma_2^+ = \emptyset$.
Then $\Pr[E_\sigma] \leq \Pr[\Q_{\tau \nu(n+1)}^F(B_\sigma^+) = 0]$
and since the upper half of $B_\sigma^+ $ is contained in
$[0,1]^2$, in this case 
\bea
\Pr[E_\sigma] \leq \exp( - \tau \nu(n) \pi r_{\nu(n)}^2 [
 ((1- 5 \eps)^2/4 ) - \eps/2] )
\nonumber \\
\leq \nu(n)^{-\alpha \tau (1-12 \eps)/4}
\leq \nu(n)^{-(1+\eps)/4}.
\lbl{0903d}
\eea
Now suppose $\sigma_2^+ \neq \emptyset$.
  Let $q_\ell$ be the last element (in the lexicographic
order) of  $\sigma_2^+$. Then
\bean
\Pr[E_\sigma] \leq \Pr[\Po_{\nu(n+1)}^F(\{q_{\ell} \}
 \oplus B_+ ( r_{\nu(n)}(1- 3 \eps) ) )  = 0] \\
\leq \exp ( - {\nu(n)} \pi r_{\nu(n)}^2 (1- 3 \eps)^2/4 ) \leq
\nu(n)^{-\alpha (1- 6 \eps)/4} \leq
\nu(n)^{-1/20},
\eean
where for the last inequality we used the fact
that $\alpha > 1/4$ and $ \eps <1/99$.
Together with \eq{0903d} this demonstrates \eq{Ani3}.
$\qed$ \\

For $m,n \in \N$, and $r >0$,
let $\cK_{n,m}(r)$ be the class of bipartite
point sets $(\X,\Y)$ in $[0,1]^2$ such that $G(\X,\Y,r)$
has at 
least one component, the vertex-set of which
has projection onto $\LL_n$ of order $m$ and contains
at least one element of $\X$.
\begin{lemm}
\lbl{Lem0816a}
Let $m \in \N$. Then almost surely,
for all but finitely many $n \in \N$ we have
$ (\Po_{n'}^F,\Q_{\tau {n'}}^F ) \notin \cK_{\nu(n),m}(r_{n'}) $
for all $n' \in \N \cap [\nu(n),\nu(n+1))$.
\end{lemm}
{\em Proof.}
By Lemmas \ref{conlem1} and \ref{AniLem},
for $n \geq n_1$ we have
\bean
\Pr[ 
\cup_{\nu(n) \leq n' < \nu(n+1) } 
\{ (\Po_n^F,\Q_{\tau n}^F ) \in \cK_{\nu(n),m}(r_{n'}) \} ]
\leq
\sum_{\sigma \in \cA_{\nu(n),m} } \Pr(E_\sigma)
\\
\leq 
|\A_{\nu(n),m}^{2} | \times \nu(n)^{-(1+\eps)}
+
|\A_{\nu(n),m}^{1} | \times \nu(n)^{-(1+\eps)/2}
+
|\A_{\nu(n),m}^{0} | \times \nu(n)^{-1/20}
\eean
and using Lemma \ref{countlem} and the
 definition $\nu(n):= n^{\lceil 4/\eps_0 \rceil} $,
and recalling that $\eps = \eps_n \geq \eps_0$  as
described just after (\ref{epscond2}),
we find that this probability
 is $O(n^{-2})$, so is summable in $n$; then the result follows
by the Borel-Cantelli lemma. $\qed$

\begin{lemm}
\lbl{Unilem}
{\em (see \cite[Lemma 9.1]{Penbk}.)}
For any two closed connected subsets $A,B$ of $[0,1]^2$
with union $A \cup B = [0,1]^2 $, the intersection
$A \cap B$ is connected.
\end{lemm}

Given $n \in \N$, let $k(n) $ be the choice of
$k \in \N$ satisfying
 $\nu(k) \leq n < \nu(k+1)$.
Also, given $K\in \N$, let $F_K(n)$
be the event that $G(\Po_n^F,\Q_{\tau n}^F,r_n)$ has two or more components
with projections onto $\LL_{\nu(k(n))}$ of  order greater than
$K$.
\begin{lemm}
\lbl{Lem0816b}
There exists $K \in \N$ such that with probability 1
the event $F_K(n) $ occurs for only finitely many $n$.
\end{lemm}
{\em Proof.}
Suppose $F_K(n)$ occurs. Then there exist distinct components
 $U=(U_1,U_2) $,
and $V = (V_1,V_2) $ in $G(\Po_n^F,\Q^F_{\tau n},r_n)$, 
both with projections onto $\LL_{\nu(k(n))}$ of order greater than $K$.
Let $U'$ be the union of closed Voronoi cells  in $[0,1]^2$
(relative to $\Po_n \cup \Q_{\tau n}$) of vertices of $U$,
and let $V'$ be the union of closed Voronoi cells in $[0,1]^2$
of vertices of $V$.

The interior of $U'$ and the interior of $ V'$ are
 disjoint subsets of  $[0,1]^2$, and we now show that
they are connected sets.
Suppose $X \in U_1, Y \in U_2$
with $\|X-Y\| \leq r_n$; then we claim the entire 
line segment $[X,Y]$ is contained in the interior
of $U'$. Indeed, 
let $z \in [X,Y]$, and suppose $z$ lies in the closed
Voronoi cell of some   $W \in \Po_n^F \cup \Q_{\tau n}^F$.
If $W \in \Po_n^F$ then 
$$
\|W -Y \| \leq \| W - z \| + \|z-Y\| \leq  
\| X - z \| + \|z-Y\|  = \|X-Y \| \leq r_n
$$  
so $W \in U$. Similarly, if $W \in \Q_{\tau n}^F$ then
$\|W-X\| \leq r_n$ so 
again $W \in U$. Hence the interior of $U'$ is connected,
and likewise for $V'$.

Let $\tilde{V}$ be the closure of the component of $[0,1]^2 \setminus U'$,
containing the interior of  $V'$, and let $\tilde{U}$ be the closure of
$[0,1]^2 \setminus \tilde{V}$ (essentially, this is the
set obtained by filling in the holes of $U'$ that are not connected to $V'$).

Then $\tilde{U},\tilde{V}$ are closed connected sets, whose union
is $[0,1]^2$.
Therefore by Lemma \ref{Unilem}, the set
 $\partial U : = \tilde{U} \cap \tilde{V}$
 is connected.
Note that $\partial U$ 
is part of the boundary of $U'$
(it is the `exterior boundary' of $U'$ relative to $V'$).

Let $T $ be the set of cube centres $q_i \in \LL_{\nu(k(n))}$ such
 that $Q_i \cap 
(\partial U) \neq \emptyset$.
Then $T$ is $*$-connected in $\cL_{\nu(k(n))}$,  i.e.
for any $x,y \in  T$,
there is a path $(x_0,x_1,\ldots,x_k)$
with $x_0=x$, $x_k=y$ and $x_i \in \LL_{\nu(k(n))}$ and
 $\|x_i -x_{i-1}\|_\infty = \eps r_{\nu(k(n))}$
for $1 \leq i \leq k$ (here $\eps = \eps_{\nu(k(n))}$.) 

Also, for each $q_i \in T$
we claim $\Po_n(Q_i) \Q_{\tau n}(Q_i) = 0$.
Indeed, suppose on the contrary 
that $\Po_n(Q_i) \Q_{\tau n}(Q_i) > 0$.
Then all points of $(\Po_n \cup \Q_{\tau n}) \cap Q_i$ 
lie in the same component of $G(\Po_n^F,\Q_{\tau n}^F,r_n)$.
If they are all in $U$, then
 $Q_i$, and all
 neighbouring $Q_j$ (including diagonal neighbours)
are contained  in $U'$. 
If all points of $(\Po_n \cup \Q_{\tau n}) \cap Q_i$ 
are not in $U$, then 
$Q_i$, and 
 all neighbouring $Q_j$ (including diagonal neighbours)
are disjoint from $U'$.
Therefore
 $(\partial U ) \cap Q_i = \emptyset$.

We now prove the  isoperimetric inequality
\bea
|T| \geq   (K/2)^{1/2}.
\lbl{isoperi}
\eea 
To see this, define the {\em width} of a nonempty closed set
$A \subset [0,1]^2$ to be the maximum difference
between $x$-coordinates of points in $A$, and
 the {\em height} of $A$
 to be the maximum difference
between $y$-coordinates of points in $A$.

We claim that
either the height or the width of $\partial U$ is 
at least $(K/2)^{1/2} \eps r_{\nu(k(n))}$. Indeed,   if not, then
$\partial U$ is contained in some square of side
$(K/2)^{1/2} \eps r_{\nu(k(n))}$, and then either $U'$ or $V'$ is contained
in that square, so either $U$ or $V$ is contained in
that square, contradicting the assumption that the projections
of $U$ and of $V$ onto $\LL_{\nu(k(n))}$ have order greater than
 $K$.  
For example, if the projection of $U$ has order greater
than $K$, then at least one of $U_1$ and $U_2$, say $U_1$,
has projection of order greater than $K/2$, and then
the union of squares of side $\eps r_{\nu(k(n))}$ centred
at vertices in the projection of $U_1$ has total area
greater than $(K/2) \eps^2 r_{\nu(k(n))}^2$, so is
not contained in  any square of side $(K/2)^{1/2} \eps r_{\nu(k(n))}$. 
Thus the claim holds, and
then \eq{isoperi} follows by
 the $*$-connectivity of $T$. 

For $\nu,m \in \N$,
let $\A'_{\nu,m}$ be the set of $*$-connected subsets of $\cL_\nu$
with $m$ elements.
By a similar argument to the proof of Lemma \ref{countlem}
(see also \cite[Lemma 9.3]{Penbk}),
there  are finite constants $\gamma$ and $C$ such that
for all $\nu, m \in \N$ we have
\bea
|\A'_{\nu,m} | \leq C (\nu/\log \nu) \gamma^m.
\lbl{0903e}
\eea
Set 
$
\phi_n : = \Pr[\Po_n(Q_i) \Q_{\tau n}(Q_i)  = 0  ] 
$ (which does not depend on $i$). By the union bound,
and \eq{0818a}, 
\bean
\phi_n \leq \exp(- n (\eps r_{\nu(k(n))})^2) + \exp( - \tau n
 (\eps r_{\nu(k(n))})^2 )
\\
\leq 2  
 \exp[ - (\tau \wedge 1)  \eps^2 (\alpha/\pi)  
(n  \log \nu(k(n) )/\nu(k(n)) ) ] 
\\
\leq 2 \nu(k(n))^{-(\tau \wedge 1) \eps^2 \alpha /\pi}
\leq 3 n^{-(\tau \wedge 1) \eps^2 \alpha /\pi},
\eean
where the last inequality holds for all large enough $n$.
Using
 \eq{isoperi} and \eq{0903e} we obtain that
\bean
\Pr[F_K(n) ] \leq \sum_{m \geq  (K/2)^{1/2}}
C (\nu(k(n))/\log \nu(k(n))) \gamma^m \phi_n^m
\\
\leq 2 C  n ( 3 \gamma n^{-\eps^2 \alpha 
(\tau \wedge 1)/\pi})^{( K/2)^{1/2} },
\eean
which is summable in $n$ provided $K$ is chosen 
so that $\eps^2 \pi^{-1}  \alpha (\tau \wedge 1)(K/2)^{1/2} >3$.
The result then follows by the Borel-Cantelli lemma.
$\qed$ \\

{\em Proof of Theorem \ref{connth2}.}
Choose $K\in \N$ as in Lemma
\ref{Lem0816b}. Writing `i.o.' for `for infinitely many $n$' 
(i.e. infinitely often), we have
\bean
 \Pr [ \G^1 (n,\tau,r_n) \notin \cK \mbox{ i.o.}]  \leq 
 \left( \sum_{m=1}^K \Pr[(\Po_n^F,\Q_{\tau n}^F) \in \cK_{\nu(k(n)),m} (r_n)
\mbox{ i.o.} ] \right) 
+ \Pr[F_K(n) \mbox{ i.o.}].
\eean
By Lemmas \ref{Lem0816a} and \ref{Lem0816b}, this is zero. $\qed$

\end{document}